\documentclass[twoside,leqno,10pt, A4]{amsart}
\usepackage{amsfonts}
\usepackage{amsmath}
\usepackage{amscd}
\usepackage{amssymb}
\usepackage{amsthm}
\usepackage{amsrefs}
\usepackage{latexsym}
\usepackage{mathrsfs}
\usepackage{bbm}
\usepackage{color}
\begin{document}

\newtheorem{theorem}[subsection]{Theorem}
\newtheorem{proposition}[subsection]{Proposition}
\newtheorem{lemma}[subsection]{Lemma}
\newtheorem{corollary}[subsection]{Corollary}
\newtheorem{conjecture}[subsection]{Conjecture}
\newtheorem{prop}[subsection]{Proposition}
\numberwithin{equation}{section}
\newcommand{\mr}{\ensuremath{\mathbb R}}
\newcommand{\mc}{\ensuremath{\mathbb C}}
\newcommand{\dif}{\mathrm{d}}
\newcommand{\intz}{\mathbb{Z}}
\newcommand{\ratq}{\mathbb{Q}}
\newcommand{\natn}{\mathbb{N}}
\newcommand{\comc}{\mathbb{C}}
\newcommand{\rear}{\mathbb{R}}
\newcommand{\prip}{\mathbb{P}}
\newcommand{\uph}{\mathbb{H}}
\newcommand{\fief}{\mathbb{F}}
\newcommand{\majorarc}{\mathfrak{M}}
\newcommand{\minorarc}{\mathfrak{m}}
\newcommand{\sings}{\mathfrak{S}}
\newcommand{\fA}{\ensuremath{\mathfrak A}}
\newcommand{\mn}{\ensuremath{\mathbb N}}
\newcommand{\mq}{\ensuremath{\mathbb Q}}
\newcommand{\half}{\tfrac{1}{2}}
\newcommand{\f}{f\times \chi}
\newcommand{\summ}{\mathop{{\sum}^{\star}}}
\newcommand{\chiq}{\chi \bmod q}
\newcommand{\chidb}{\chi \bmod db}
\newcommand{\chid}{\chi \bmod d}
\newcommand{\sym}{\text{sym}^2}
\newcommand{\hhalf}{\tfrac{1}{2}}
\newcommand{\sumstar}{\sideset{}{^*}\sum}
\newcommand{\sumprime}{\sideset{}{'}\sum}
\newcommand{\sumprimeprime}{\sideset{}{''}\sum}
\newcommand{\shortmod}{\ensuremath{\negthickspace \negthickspace \negthickspace \pmod}}
\newcommand{\V}{V\left(\frac{nm}{q^2}\right)}
\newcommand{\sumi}{\mathop{{\sum}^{\dagger}}}
\newcommand{\mz}{\ensuremath{\mathbb Z}}
\newcommand{\leg}[2]{\left(\frac{#1}{#2}\right)}
\newcommand{\muK}{\mu_{\omega}}

\title[Moments {H}ecke {$L$}-functions with quadratic characters]{Moments of Quadratic {H}ecke {$L$}-functions of Imaginary Quadratic Number Fields}

\date{\today}
\author{Peng Gao and Liangyi Zhao}

\begin{abstract}
In this paper, we study the moments of central values of Hecke $L$-functions associated with quadratic characters in $\mq(i)$ and $\mq(\omega)$ with $\omega=\exp(2\pi i/3)$ and establish some quantitative non-vanishing result for these $L$-values.
\end{abstract}

\maketitle

\noindent {\bf Mathematics Subject Classification (2010)}: 11M41, 11L40  \newline

\noindent {\bf Keywords}: quadratic Hecke characters, Hecke $L$-functions

\section{Introduction}

The values of $L$-functions at the central point of its symmetry encode a lot of arithmetic information.  For example, the Birch-Swinnerton-Dyer conjecture asserts that the algebraic rank of an elliptic curve equals the order of vanishing of the $L$-function associated with the curve at its central point.  The average value of $L$-functions at $s=1/2$ over a family of characters of a fixed order has been a subject of extensive study.  M. Jutila \cite{Jutila} gave the evaluation of the mean value of $L(1/2, \chi)$ for quadratic Dirichlet characters.  The error term in the asymptotic formula in \cite{Jutila} was later improved in \cites{DoHo, MPY, ViTa}.  S. Baier and M. P. Young \cite{B&Y} studied the moments of $L(1/2, \chi)$ for cubic Dirichlet characters.  Literature also abounds in the investigation of moments of Hecke $L$-functions associated with various families of characters of a fixed order \cites{FaHL, GHP, FHL, Luo, Diac, G&Zhao1}.   In this paper, we shall investigate the first and second moments of the central values of a family of $L$-functions associated with quadratic Hecke characters.  From these results, we deduce results regarding the non-vanishing of these $L$-functions at $s=1/2$. \newline

Set $\omega = \exp (2\pi i/3)$.  Let $K =\mq(i)$ or $\mq(\omega)$ and $\chi$ be a primitive Hecke character $\pmod {m}$ of trivial infinite type defined on $\mathcal{O}_K$.  The Hecke $L$-function associated
with $\chi$ is defined for $\Re(s) > 1$ by
\begin{equation*}
  L(s, \chi) = \sum_{0 \neq \mathcal{A} \subset
  \mathcal{O}_K}\chi_c(\mathcal{A})(N(\mathcal{A}))^{-s},
\end{equation*}
  where $\mathcal{A}$ runs over all non-zero integral ideals in $K$ and $N(\mathcal{A})$ is the
norm of $\mathcal{A}$. As shown by E. Hecke, $L(s, \chi)$ admits
analytic continuation to an entire function and satisfies the
functional equation (\cite[Theorem 3.8]{HIEK})
\begin{align} \label{1.1}
  \Lambda(s, \chi) = W(\chi)(N(m))^{-1/2}\Lambda(1-s, \overline{\chi}),
\end{align}
   where $|W(\chi)|=(N(m))^{1/2}$ and
\begin{align*}
  \Lambda(s, \chi) = (|D_K|N(m))^{s/2}(2\pi)^{-s}\Gamma(s)L(s, \chi),
\end{align*}
   with $D_K$ being the discriminant of $K$.  So $D_K=-4$ if $K= \mq(i)$ and $D_K=-3$ if $K=\mq(\omega)$. \newline

   For a square-free, non-unit $c \in \intz[i]$ congruent to $1 \pmod {16}$, let $\chi_c=(\frac {\cdot}{c})_4$ be the quartic
residue symbol defined in Section \ref{sec2.4}. It gives rise to a primitive Hecke character $\pmod {c}$ of trivial infinite type.
In this case, the functional equation \eqref{1.1} holds with
\begin{align*}
%%\label{1.2}
   W(\chi_c)=\sum_{a \in \mathcal{O}_{K}/(c)}\chi_c(a)e\Big ( \text{Tr}\Big (\frac {a}{\delta c}\Big )\Big ),
\end{align*}
   with $(\delta)=(\sqrt{-4})$ being the different of $K$ and $e(z) = \exp( 2 \pi i z)$. Note that
   $W(\chi_c)$ equals the Gauss sum  $g(c)$ defined in
   Section \ref{sec2.4}. \newline

In \cite{G&Zhao1}, we studied the first moment of  $L(1/2,\chi_c)$.  Our approach was to decompose $L(1/2, \chi_c)$ using the approximate functional equation into two sums.  The treatment on one of them relied crucially on an estimate of a smooth Gauss sum. We then applied a result of S. J. Patterson  \cite[Lemma, p. 200]{P} to show that the contribution of the corresponding term was small, yielding an acceptable error term. \newline

   The result of Patterson is more general. In fact, for any algebraic number field $F$, let $\mu_n(F)$ be the set of $n$-th roots of unity in $F$ with $n>2$. Suppose that the cardinality of $\mu_n(F)$ is $n$. For any $n$-th order Hecke character $\chi$ and any injective character $\varepsilon: \mu_n(F) \rightarrow \mc^{\times}$, one can define a Gauss sum for $\varepsilon(\chi)$ in a manner similar to the Gauss sum defined in
   Section \ref{sec2.4}.  Patterson's result then gives that these Gauss sums have a lot of cancellation on average. On the other hand, one does not expect such a result to hold when $\varepsilon$ is not injective. For example, when $F =\mq(i)$, let $\varepsilon$ be the character: $x \mapsto x^2$ from $\mu_4(K)$ to $\mc^{\times}$. Then for any quartic symbol $\chi$, $\varepsilon(\chi)$ becomes a quadratic symbol defined in Section \ref{sec2.4} and the associated Gauss sum becomes a quadratic Gauss sum. To fathom the behavior of these quadratic Gauss sums on average, one can examine their rational analogues: the Gauss sums for quadratic Dirichlet characters. In the work of M. Jutila \cite{Jutila} and K. Soundararajan \cite{sound1} on the first moment of  $L(1/2,\chi)$ for quadratic Dirichlet characters $\chi$, it is shown that both sums coming from the approximate functional equation contribute to the main term, unlike the case of quartic symbols discussed above. \newline

   Motivated by these observations, it is our goal in this paper to study the first and second moments of Hecke $L$-functions with quadratic characters in $\mq(i)$ and $\mq(\omega)$ at the central point. In this case Patterson's result is no longer in display and our approach is similar to that used in \cite{sound1}. In particular, we shall see that both sums coming from the approximate functional equation contribute to the main term, just as the case for the quadratic Dirichlet $L$-functions. \newline

  Let $\Phi$ be a smooth Schwarz class function compactly supported in $(1, 2)$ and we assume that $0 \leq \Phi(t) \leq 1$ for all $t$. We have
\begin{theorem}
\label{firstmoment}
Let $K = \ratq(i)$ or $\ratq(\omega)$.  For $y \rightarrow \infty$ and any $\varepsilon > 0$, we have
\begin{align} \label{1stmom}
   \sumstar_{c \in \mathcal{O}_K } L\left( \frac{1}{2},
   \chi_c \right)\Phi\left( \frac{N(c)}{y} \right)=A_K B_K \frac{\pi^2}{\zeta_{K}(2)} \hat{\Phi}(0) y\log y +C_K \hat{\Phi}(0)y+O\left( y^{(3+\theta)/4} \right),
\end{align}
   where $\zeta_{K}(s)$ is the Dedekind zeta function of $K$, $\theta=131/416$,
\begin{align} \label{1.4}
   A_{\ratq(i)}= \prod_{\substack {\pi \text{ prime in $\mq(i)$} \\(\pi,2)=1} } \left( 1-\frac {1}{(N(\pi)+1)N(\pi)} \right) , \quad A_{\ratq(\omega)}= \prod_{\substack {\pi \text{ prime in $\mq(\omega)$} \\(\pi,6)=1} } \left( 1-\frac {1}{(N(\pi)+1)N(\pi)} \right),
\end{align}
\[ B_{\ratq(i)} = \frac{2+\sqrt{2}}{3072},  \quad B_{\ratq(\omega)} = \frac{3+\sqrt{3}}{19440} , \quad \hat{\Phi}(0) = \int\limits_1^2 \Phi(x) \dif x , \]
$C_K$ is a constant and $\sum^{*}$ indicates that the sum runs over square-free elements of $\mz[i]$ congruent to $1 \pmod {16}$ if $K = \ratq(i)$ and square-free elements of $\mz[\omega]$ congruent to $1 \pmod{36}$ if $K = \ratq(\omega)$.
\end{theorem}

We note here that $\theta$ arises from an application of the currently best known formula in the Gauss circle problem \cite{Huxley1}.  Therefore, any improvement in the study of that problem will also lead to an improvement in the $O$-term in \eqref{1stmom}. \newline

 As for the second moment, we prove
\begin{theorem} \label{secmom}
Using the same notation as Theorem~\ref{firstmoment},  for $y \rightarrow \infty$ and any $\varepsilon > 0$, we have
\begin{equation*} \label{2ndmom}
   \sumstar_{c \in \mathcal{O}_K} \left| L \left( \frac{1}{2},
   \chi_c \right) \right|^2\Phi \left(  \frac{N(c)}{y} \right) \ll_{\varepsilon} y^{1+\varepsilon}.
\end{equation*}
\end{theorem}

S. Chowla \cite{chow} was the first to conjecture that a Dirichlet $L$-function is never zero at $s=1/2$ and it is believed that an $L$-function cannot vanish at its central point unless there is a compelling reason (an elliptic curve with a postive algebraic rank or the root number being $-1$) that there should be a zero there.  In this vein, we give the following non-vanishing results for the $L$-functions under our consideration.

\begin{corollary}
  For $y \rightarrow \infty$ and any $\varepsilon > 0$, we have
\[    \# \left\{ c \in \mz[i] : c \equiv 1 \pmod {16}, \; N(c) \leq y, \; L\left( \frac{1}{2}, \chi_c \right) \neq 0 \right\} \gg_{\varepsilon}
   y^{1-\varepsilon} \]
and
\[    \# \left\{ c \in \mz[\omega] : c \equiv 1 \pmod {36}, \; N(c) \leq y, \; L\left( \frac{1}{2}, \chi_c \right) \neq 0 \right\} \gg_{\varepsilon}
   y^{1-\varepsilon} . \]
\end{corollary}

\begin{proof}
The lower bounds follow from Theorems~\ref{firstmoment} and \ref{secmom}, via standard arguments using Cauchy's inequality (see \cite{Luo}).\end{proof}

\section{Preliminaries}
\label{sec 2}

In this section, we write down the preliminary results required in the proof of our main theorems.

%%----------------------------------------------------------------------------
\subsection{Quadratic symbol and quadratic Gauss sum in $\mq(i)$}
\label{sec2.4}
%%----------------------------------------------------------------------------
   The symbol $(\frac{\cdot}{n})_4$ is the quartic
residue symbol in the ring $\mz[i]$.  For a prime $\pi \in \mz[i]$
with $N(\pi) \neq 2$, the quartic character is defined for $a \in
\mz[i]$, $(a, \pi)=1$ by $\leg{a}{\pi}_4 \equiv
a^{(N(\pi)-1)/4} \pmod{\pi}$, with $\leg{a}{\pi}_4 \in \{
\pm 1, \pm i \}$. When $\pi | a$, we define
$\leg{a}{\pi}_4 =0$.  Then the quartic character can be extended
to any composite $n$ with $(N(n), 2)=1$ multiplicatively. We further define $(\frac{\cdot}{n})=\leg {\cdot}{n}^2_4$  to be the quadratic
residue symbol for all $n \in \mz[i]$ with $(N(n), 2)=1$. \newline

 Note that in $\intz[i]$, every ideal coprime to $2$ has a unique
generator congruent to 1 modulo $(1+i)^3$.  Such a generator is
called primary. Observe that a non-unit
$n=a+bi$ in $\mz[i]$ is congruent to $1
\bmod{(1+i)^3}$ if and only if $a \equiv 1 \pmod{4}, b \equiv
0 \pmod{4}$ or $a \equiv 3 \pmod{4}, b \equiv 2 \pmod{4}$ by \cite[Lemma 6, p. 121]{I&R}. \newline

  Recall that (see \cite[Theorem 2, p. 123]{I&R}) the quartic reciprocity law states
that for two primary primes  $m, n \in \mz[i]$,
\begin{align*}
 \leg{m}{n}_4 = \leg{n}{m}_4(-1)^{((N(n)-1)/4)((N(m)-1)/4)}.
\end{align*}
Thus the following quadratic reciprocity law holds for two primary primes  $m, n \in \mz[i]$:
\begin{align*}
 \leg{m}{n} = \leg{n}{m}.
\end{align*}

 For a non-unit $n \in \mz[i]$, the quadratic Gauss sum $g(n)$ is defined by
\[    g(n) =\sum_{x \bmod{n}} \leg{x}{n} \widetilde{e}\leg{x}{n}, \; \; \; \mbox{where} \; \; \; \widetilde{e}(z) =\exp \left( 2\pi i  \left( \frac {z}{2i} - \frac {\overline{z}}{2i} \right) \right). \]

     From the supplement theorem to the quartic reciprocity law (see for example, Lemma 8.2.1 and Theorem 8.2.4 in \cite{BEW}),
we have for $n=a+bi$ primary,
\begin{align*}
%%\label{2.05}
  \leg {i}{n}_4=i^{(1-a)/2} \qquad \mbox{and} \qquad  \hspace{0.1in} \leg {1+i}{n}_4=i^{(a-b-1-b^2)/4}.
\end{align*}
   It follows that for any $c \equiv 1 \pmod {16}$, we have
\begin{align}
\label{2.5i}
  \leg {i}{c}_4=\leg {1+i}{c}_4=1.
\end{align}

   This shows that $\chi_c =(\frac {\cdot}{c})_4$ is trivial on
units, hence for any $c$ square-free and congruent to $1 \pmod {16}$, $\chi_c$ can be regarded as a primitive character of the ray
class group $h_{(c)}$. We recall here that for any $c$, the ray
class group $h_{(c)}$ is defined to be $I_{(c)}/P_{(c)}$, where
$I_{(c)} = \{ \mathcal{A} \in I : (\mathcal{A}, (c)) = 1 \}$ and
$P_{(c)} = \{(a) \in P : a \equiv 1 \pmod{c} \}$ with $I$ and $P$
denoting the group of fractional ideals in $K$ and the subgroup of
principal ideals, respectively. \newline

%\subsection{Evaluation of the Gauss sums}
%\label{section:Gauss}

We shall determine the exact value of $g(c)$ where $c \equiv 1 \pmod {16}$ is square-free in $\mz[i]$. We have $g(1)=1$ by definition. If $c \neq 1$, we can write $c=\pi_1\cdots \pi_k$ with $\pi_i \equiv 1 \pmod {(1+i)^3}$ being distinct primes. It follows from quadratic reciprocity that
\begin{align*}
  g(c)=  \prod^k_{i=1}g(\pi_i).
\end{align*}

   Thus, it suffices to compute $g(\pi)$ for a prime $\pi \equiv 1 \pmod {(1+i)^3}$. This evaluation is available in \cite[Proposition 2.2]{Onodera} (be aware that the definition of the Gauss sum in \cite{Onodera} is different from the one here) and we have
\[   g(\pi)=\leg{-1}{\pi}_4N(\pi)^{1/2}. \]

We conclude readily from the above discussions and \eqref{2.5i} that for $c \equiv 1 \pmod {16}$ and square-free,
\begin{equation} \label{gcvalueQi}
   g(c)=\leg {-1}{c}_4N(c)^{1/2}=N(c)^{1/2}.
\end{equation}

\subsection{Quadratic symbol, Kronecker symbol and quadratic Gauss sum in $\mq(\omega)$}
   It is well-known that $K=\mq(\omega)$ has class number $1$, and the symbol $(\frac{\cdot}{n})$ is the quadratic
residue symbol in the ring of integers $\mathcal{O}_K =\mz[\omega]$.  For a prime $\pi \in \mz[\omega], \pi \neq 2$, the quadratic character is defined for $a \in
\mz[\omega]$, $(a, \pi)=1$ by $\leg{a}{\pi} \equiv
a^{(N(\pi)-1)/2} \pmod{\pi}$, with $\leg{a}{\pi} \in \{
\pm 1 \}$. When $\pi | a$, we define
$\leg{a}{\pi} =0$.  Then the quadratic character can be extended
to any composite $n$ with $(N(n), 2)=1$ multiplicatively. We further define $\leg {\cdot}{n}=1$ when $n$ is a unit in $\mz[\omega]$. \newline

In $\intz[\omega]$, every ideal co-prime to $3$ has a unique generator congruent to 1  $\pmod 3$.  Such a generator is
called primary.  Observe that a non-unit
$n=a+b\omega$ in $\mz[\omega]$ is congruent to $1
\pmod{3}$ if and only if $a \equiv 1 \pmod{3}$, and $b \equiv
0 \pmod{3}$ (see the discussions before \cite[Proposition 9.3.5]{I&R}). \newline

  We shall say that any $n \in \mz[\omega]$ is cubic $E$-primary if $n^3=a+b\omega$ with $a, b \in \mz$ such that $6 | b$ and $a+b \equiv 1 \pmod 4$. Any cubic $E$-primary number is thus co-prime to $2$. \newline

    It follows from \cite[Lemma 7.9]{Lemmermeyer} that any $n=a+b\omega \in \mz[\omega]$ is cubic $E$-primary if and only if
\begin{align}
\label{cubicE}
   & a+b \equiv 1 \pmod 4, \quad  \text{if} \quad 2 | b,  \\
   & b \equiv 1 \pmod 4, \quad  \text{if} \quad 2 | a,  \nonumber \\
   & a \equiv 3 \pmod 4, \quad  \text{if} \quad 2 \nmid ab.  \nonumber
\end{align}

     Furthermore, the following quadratic reciprocity law holds for two cubic $E$-primary, co-prime numbers $n, m \in \mz[\omega]$ :
\begin{align}
\label{quadreciQw}
    \leg {n}{m} =\leg{m}{n}(-1)^{((N(n)-1)/2)((N(m)-1)/2)}.
\end{align}

  Let $c \in \mz[\omega]$, we say that $c$ is $E$-primary if we can write it as $c=(-(1-\omega))^rc'$ with $r \geq 0, r \in \mz, (c', 6)=1$, $c'$ is cubic $E$-primary and either $c'$ or $-c'$ is primary. Note that our definition of $E$-primary here is the same as that defined in \cite[Section 7.3]{Lemmermeyer} when $(c,6)=1$. One checks easily that every ideal co-prime to $2$ in $\intz[\omega]$ has a unique $E$-primary generator. \newline

  One also has the following supplementary laws for $n=a+b\omega, (n,6)=1$ being $E$-primary (see \cite[Theorem 7.10]{Lemmermeyer}),
\begin{align}
\label{2.05}
  \leg {-1}{n}=(-1)^{(N(n)-1)/2}, \qquad \leg {1-\omega}{n}=\leg {a}{3}_{\mz} \qquad \mbox{and} \qquad  \hspace{0.1in} \leg {2}{n}=\leg {2}{N(n)}_{\mz},
\end{align}
   where $\leg {\cdot}{\cdot}_{\mz}$ denotes the Jacobi symbol in $\mz$. One checks that the last expression above holds in fact for all $E$-primary numbers. \newline

When $c \in \mz[\omega]$ which is square-free and congruent to $1 \pmod {36}$, it follows from \eqref{2.05} that
\begin{align}
\label{2.5}
  \leg {-1}{c}=\leg {1-\omega}{c}=\leg {2}{c}=1.
\end{align}
This shows that $\chi_c =(\frac {\cdot}{c})$ is trivial on units and it can be regarded as a primitive character of the ray class group $h_{(c)}$. \newline

 For any element $c \in \mz[\omega], (c, 2)=1$, we can define a quadratic Dirichlet character $\chi^{(-8c)} \pmod {8c}$ such that for any $n \in (\mz[\omega]/(8c\mz[\omega]))^*$,
\begin{align*}
   \chi^{(-8c)}(n)=\leg {-8c}{n}.
\end{align*}

    One deduces from \eqref{2.05} and the quadratic reciprocity that  $\chi^{(-8c)}(n)=1$ when $n \equiv 1 \pmod {8c}$. It follows from this that  $\chi^{(-8c)}(n)$ is well-defined. As $\chi^{(-8c)}(n)$ is clearly multiplicative and of order $2$ and is trivial on units, it can be regarded as a quadratic Hecke character $\pmod {8c}$ of trivial infinite type. We denote $\chi^{(-8c)}$ for this Hecke character as well and we call it the Kronecker symbol. Furthermore, when $c$ is square-free, $\chi^{(-8c)}$ is non-principal and primitive. To see this, we write $c=u_c \cdot \varpi_1 \cdots \varpi_k$ with $u_c$ a unit and $\varpi_j$  being $E$-primary primes. Suppose $\chi^{(-8c)}$ is induced by some $\chi$ modulo $c'$ with $\varpi_j \nmid c'$, then by the Chinese Remainder Theorem, there exists an $n$ such that $n \equiv 1 \pmod {8c/\varpi_j}$ and $\leg {\varpi_j}{n} \neq 1$. It follows that $\chi(n)=1$ but $\chi^{(-8c)}(n) \neq 1$, a contradiction. Thus, $\chi^{(-8c)}$ can only be possibly induced by some $\chi$ modulo $4c$. By the Chinese Remainder Theorem, there exists an $n$ such that $n \equiv 1 \pmod {c}$ and $n \equiv 1+4\omega \pmod {8}$. As this $n \equiv 1 \pmod {4}$, it follows that $n \equiv 1 \pmod {4c}$, hence $\chi(n)=1$ but $\chi^{(-8c)}(n)=\leg {2}{n}=-1 \neq 1$ (note that $\leg {u}{n}=1$ when $u$ is a unit in $\mz[\omega]$) and this implies that $\chi^{(-8c)}$ is primitive. This also shows that $\chi^{(-8c)}$ is non-principal. \newline

For a non-unit $n \in \mz[\omega]$, $(n,2)=1$, the quadratic Gauss sum $g(n)$ is defined by
\[    g(n) =\sum_{x \bmod{n}} \leg{x}{n} \widetilde{e}\leg{x}{n} , \quad \mbox{where} \quad \widetilde{e}(z) =e \left( \frac {z}{\sqrt{-3}} - \frac {\overline{z}}{\sqrt{-3}} \right). \]

   It follows from the definition that $g(1)=1$. The following well-known relation (see \cite{Diac}) now holds for all $n$:
\begin{align*}
%%\label{2.1}
   |g(n)|& =\begin{cases}
    \sqrt{N(n)} \qquad & \text{if $n$ is square-free}, \\
     0 \qquad & \text{otherwise}.
    \end{cases}
\end{align*}

   The following properties of $g(n)$ can be easily derived from the definition:
\begin{align}
\label{2.03}
   g(n_1 n_2) =\leg{n_2}{n_1}\leg{n_1}{n_2}g(n_1) g(n_2), \quad (n_1, n_2) = 1.
\end{align}

   In what follows we compute the value of $g(c)$ where $c \equiv 1 \pmod {36}$ is square-free in $\mz[\omega]$.
We first evaluate the Gauss sum at each $E$-primary prime $\varpi$. We have the following
\begin{lemma}
\label{Gausssum}
   Let $\varpi$ be an $E$-primary prime in $\mz[\omega]$. Then
\begin{align*}
   g(\varpi)= \begin{cases}
    N(\varpi)^{1/2} \qquad & \text{if} \qquad N(\varpi) \equiv 1 \pmod 4,\\
    -i N(\varpi)^{1/2} \qquad & \text{if} \qquad N(\varpi) \equiv -1 \pmod 4.
\end{cases}
\end{align*}
\end{lemma}
\begin{proof}
   Note that $-(1-\omega)$ is an $E$-primary prime lying above the rational prime $3$ and direct computation shows that $g(-(1-\omega))=-\sqrt{3}i$. Now we consider the case when $N(\varpi)=p$ is a prime $\equiv 1 \pmod 3$ in $\mz$.  It follows that $\overline{\varpi}$ is also a prime in $\mz[\omega]$ such that $(\varpi, \overline{\varpi}) = 1$.  Here we use $\bar{n}$ to denote the
complex conjugate of $n$, for any $n \in \mc$. We consider the Gauss sum associated to the quadratic Dirichlet character $\chi_p=\leg {\cdot}{p}_{\mz}$:
\begin{equation*}
%%\label{2.08}
 \tau(\chi_p) =  \sum_{1 \leq x \leq p} \left(\frac{x}{p}\right)_{\mz} e \left( \frac{x}{p} \right) =  \sum_{1 \leq x \leq N(\varpi)} \left(\frac{x}{\varpi}\right) e \left( \frac{x}{N(\varpi)} \right).
\end{equation*}

  Now write $x = y \overline{\varpi} + \overline{y} \varpi$, where $y$ varies over a set
of representatives in $\mz[\omega] \pmod{\varpi}$, then it is easy to see that as $y$ varies
$\pmod{\varpi}$, $x$ varies $\pmod{N(\varpi)}$ in $\mz$.  We deduce that
\begin{equation*}
%%\label{2.09}
  \tau(\chi_p) =  \sum_{y \bmod{\varpi}} \left(\frac{y \overline{\varpi}}{\varpi}\right) e \left( \frac{y}{\varpi} +\overline{\frac {y}{\varpi} } \right)  =
 \left(\frac{\overline{\varpi}}{\varpi}\right)\sum_{y \bmod{\varpi}} \left(\frac{y}{\varpi}\right) e \left( \frac{y}{\varpi} +\overline{\frac {y}{\varpi} } \right).
\end{equation*}

Observe that
\begin{equation*}
%%\label{2.10}
\begin{split}
  g(\varpi) & =  \sum_{x \bmod{\varpi}} \leg{x}{\varpi} e \left( \frac{1}{\sqrt{-3}}  \left( \frac {x}{\varpi} - \frac{\overline{x}}{\varpi} \right) \right)= \leg{\sqrt{-3}}{\varpi}\sum_{x \bmod{\varpi}} \leg{x}{\varpi} e \left( \frac {x}{\varpi } + \frac {\overline{x}}{\varpi} \right)\\
  & =\left(\frac{\sqrt{-3}\overline{\varpi}}{\varpi}\right)\tau(\chi_p)=\left(\frac{\omega(1-\omega)\overline{\varpi}}{\varpi}\right)\tau(\chi_p).
  \end{split}
\end{equation*}

Next, we shall show that
\begin{equation*}
    \left(\frac{\omega(1-\omega)\overline{\varpi}}{\varpi}\right)=\left(\frac{\omega(1-\omega)\overline{\varpi}^3}{\varpi^3}\right).
\end{equation*}

  For $\varpi$ being $E$-primary, our discussion above implies that we can write $\varpi^3=a+b\omega$ with $6|b, a+b \equiv 1 \pmod 4$ and it follows that $\overline{\varpi}^3=a+b\omega^2$. Hence
\begin{equation*}
%%\label{2.011}
 \left(\frac{\omega(1-\omega)\overline{\varpi}^3}{\varpi^3}\right)=\left(\frac{(1-\omega)(b+a\omega)}{a+b\omega}\right)
 =\left(\frac{(1-\omega)(b+a\omega-a-b\omega)}{a+b\omega}\right)=\left(\frac{(1-\omega)^2(b-a)}{a+b\omega}\right)=\left(\frac{b-a}{a+b\omega}\right).
\end{equation*}

    Note that when $a+b\omega$ satisfies $6|b, a+b \equiv 1 \pmod 4$, $a-b$ is also $E$-primary.  Thus, it follows from \eqref{quadreciQw} that
\begin{align*}
 \left(\frac{b-a}{a+b\omega}\right)=\left(\frac{-1}{a+b\omega}\right)\left(\frac{a-b}{a+b\omega}\right)
 =\left(\frac{-1}{a+b\omega}\right)(-1)^{((N(a-b)-1)/2)((N(a+b\omega)-1)/2)}
 \left(\frac {a+b\omega}{a-b}\right).
\end{align*}

   Note that $N(a-b)=(a-b)^2 \equiv 1 \pmod 4$. We then conclude that
\begin{align*}
 \left(\frac{b-a}{a+b\omega}\right)=\left(\frac{-1}{a+b\omega}\right)
 \left(\frac {a+b\omega}{a-b}\right)=\left(\frac{-1}{a+b\omega}\right)
 \left(\frac {a+b\omega+(a-b)\omega}{a-b}\right)=\left(\frac{-1}{a+b\omega}\right)
 \left(\frac {a}{a-b}\right)\left(\frac {1+\omega}{a-b}\right).
\end{align*}

   Using the relation $1+\omega+\omega^2=0$, we see that
\begin{align*}
 \left(\frac{b-a}{a+b\omega}\right)=\left(\frac{-1}{a+b\omega}\right)
 \left(\frac {-a}{a-b}\right)\left(\frac {\omega^2}{a-b}\right)=\left(\frac{-1}{a+b\omega}\right)
 \left(\frac {-a}{a-b}\right).
\end{align*}

    We note that for two co-prime $a, b \in \mz$ (see \cite[p. 219]{Lemmermeyer}), we have
\begin{align*}
   \leg {a}{b}=1.
\end{align*}

   We then conclude from the above discussions that we have
\begin{align*}
   g(\varpi)=\left(\frac{-1}{a+b\omega}\right)\tau(\chi_p)=\leg {-1}{\varpi}\tau(\chi_p)=(-1)^{(N(\varpi)-1)/2}\tau(\chi_p).
\end{align*}
  As it follows from \cite[Chap. 2]{Da} that $\tau(\chi_p)=p^{1/2}=N(\varpi)^{1/2}$ when $p \equiv 1 \pmod 4$ and $\tau(\chi_p)=ip^{1/2}=iN(\varpi)^{1/2}$ when $p \equiv 3 \pmod 4$, this completes the proof for the case when $N(\varpi)$ is a rational prime $\equiv 1 \pmod 3$. \newline

   Next, let $p \equiv 2 \pmod 3, p \neq 2$ be a prime in $\mz$, then $p$ is also a prime in $\mz[\omega]$, we now compute $g(p)$. As in \cite[Chap. 2]{Da} (note that in this case, we still have $\sum_{x \bmod p}\tilde{e}(x/p)=0$), we have
\begin{equation*}
  g(p) =  \sum_{x \bmod{p}}\tilde{e} \left( \frac {x^2}p \right).
\end{equation*}
   We now write $x=a+b\omega$ with $a, b \pmod p$ in $\mz$ to see that
\begin{equation*}
  g(p) =  \sum_{x \bmod{p}}\tilde{e} \left( \frac {x^2}p \right)=\sum^p_{b=1}\sum^p_{a=1}e \left( \frac {2ab-b^2}{p} \right)=p=N(p)^{1/2}.
\end{equation*}
  As $N(p)=p^2 \equiv 1 \pmod 4$, this completes the proof of the lemma.
\end{proof}

Now for a fixed square-free $c \neq 1, c \equiv 1 \pmod {36}$, $c$ is both primary and $E$-primary. By writing $c$ as products of primary primes and adjusting by a possible factor of $-1$, we can write $c=\varpi_1\cdots \varpi_k$ or $-\varpi_1\cdots \varpi_k$ with $\varpi_i$ being distinct $E$-primary primes. Then $c^3=\varpi^3_1\cdots \varpi^3_k$ or $-\varpi^3_1\cdots \varpi^3_k$. We write $c^3=a+b\omega$ and note that $a,b$ satisfy \eqref{cubicE} by definition.  The same consideration for $(\varpi_1\cdots \varpi_k)^3$ enables us to conclude that $c=\varpi_1\cdots \varpi_k$ as $c$ and $c^3$ differ only by a possible factor of $-1$. As $N(c) \equiv 1 \pmod 4$, we conclude that there must be an even number of $\varpi_j$ in the decomposition of $c$ such that $N(\varpi_j) \equiv -1 \pmod 4$. We may assume that $\varpi_1, \cdots, \varpi_{2k_0}$ are such primes. It follows from \eqref{2.03} and \eqref{quadreciQw} and Lemma \ref{Gausssum} that
\begin{align*}
  g(c)=  g(\varpi_1\cdots \varpi_{2k_0}) \prod^k_{j=2k_0+1}g(\varpi_j)=g(\varpi_1\cdots \varpi_{2k_0}) N \left( \prod^k_{j=2k_0+1}\varpi_j \right)^{1/2}.
\end{align*}

   Using Lemma \ref{Gausssum} and induction on $k_0$ shows that
\[  g(\varpi_1\cdots \varpi_{2k_0})=N \left( \prod^{2k_0}_{j=1}\varpi_j \right)^{1/2}. \]

   We then conclude that
\begin{equation} \label{gcvalueQomega}
  g(c)= N(c)^{1/2}.
\end{equation}

\subsection{The approximate functional equation}

   Let $\chi$ be a primitive Hecke character $\pmod {m}$ of trivial infinite type. Let $G(s)$ be any even function which is holomorphic and bounded in the strip $-4<\Re(s)<4$ satisfying $G(0)=1$. For $t \in \rear$, by evaluating the integral
\begin{equation*}
   \frac {1}{2\pi
   i}\int\limits\limits_{(2)}(2\pi)^{-(s+1/2+it)}\Gamma \left( s+\frac{1}{2}+it \right)L \left( s+\frac{1}{2}+it, \chi \right)G(s) x^s \frac{\dif s}{s}
\end{equation*}
   in two ways, we derive the following expression for $L(1/2+it, \chi)$:
\begin{equation} \label{approxfuneq}
\begin{split}
 L \left( \frac{1}{2}+it, \chi \right) = \sum_{0 \neq \mathcal{A} \subset
  \mathcal{O}_K} & \frac{\chi(\mathcal{A})}{N(\mathcal{A})^{1/2+it}}V_t \left(\frac{2\pi  N(\mathcal{A})}{x} \right) \\
  & + \frac{W(\chi)}{N(m)^{1/2}}\left(\frac {(2\pi)^2}{|D_k|N(m)} \right )^{it} \frac {\Gamma (1/2-it)}{\Gamma (1/2+it)}\sum_{0 \neq \mathcal{A} \subset
  \mathcal{O}_K}\frac{\overline{\chi}(\mathcal{A})}{N(\mathcal{A})^{1/2-it}}V_{-t}\left(\frac{2\pi
  N(\mathcal{A})x}{|D_K|N(m)} \right),
     \end{split}
\end{equation}
    where $W(\chi)$ is as in \eqref{1.1} and 
\begin{align}
\label{2.14}
  V_t \left(\xi \right)=\frac {1}{2\pi
   i}\int\limits\limits_{(2)}\frac {\Gamma(s+1/2+it)}{\Gamma (1/2+it)}G(s)\frac
   {\xi^{-s}}{s} \ \dif s.
\end{align}

    We write $V$ for $V_0$ and note that for a suitable $G(s)$ (for example $G(s)=e^{-s^2}$), we have for any $c>0$ (see \cite[Proposition 5.4]{HIEK}):
\begin{align}
\label{2.15}
  V_t \left(\xi \right) \ll \left( 1+\frac{\xi}{1+|t|} \right)^{-c}.
\end{align}

   On the other hand, when $G(s)=1$, we have (see \cite[Lemma 2.1]{sound1}) for the $j$-th derivative of $V(\xi)$,
\begin{equation} \label{2.07}
      V\left (\xi \right) = 1+O(\xi^{1/2-\epsilon}) \; \mbox{for} \; 0<\xi<1   \quad \mbox{and} \quad V^{(j)}\left (\xi \right) =O(e^{-\xi}) \; \mbox{for} \; \xi >0, j \geq 0.
\end{equation}

If $\chi_c$ is a quadratic Hecke character in $\mq(i)$, we derive readily from \eqref{approxfuneq} by setting $x=(|D_K|N(c))^{1/2}$ the following expression
\begin{equation} \label{quadapproxfuneqQi}
 L \left( \frac{1}{2}, \chi_c \right) = 2\sum_{0 \neq \mathcal{A} \subset
  \mathcal{O}_K}\frac{\chi_c(\mathcal{A})}{N(\mathcal{A})^{1/2}}V \left(\frac{\pi N(\mathcal{A}) }{N(c)^{1/2}} \right).
  \end{equation}
Note that in this situation, $W(\chi) =g(c)$ in \eqref{gcvalueQi}. \newline

A similar consideration will give that for a quadratic character $\chi_c$ in $\mq(\omega)$, 
\begin{equation} \label{quadapproxfuneqQomega}
 L \left( \frac{1}{2}, \chi_c \right) = 2\sum_{0 \neq \mathcal{A} \subset
  \mathcal{O}_K}\frac{\chi_c(\mathcal{A})}{N(\mathcal{A})^{1/2}}V \left(\frac{2 \pi N(\mathcal{A}) }{(3N(c))^{1/2}} \right).
  \end{equation}
  
The above discussions also apply to $c=1$, provided that we define $g(1)=1$ and interpret $\chi_1$ as the principal character $\pmod 1$ so that $L(s, \chi_1)$ becomes $K$, a convention we shall follow in the sequel.

\subsection{The large sieve with quadratic symbols}  The large sieve inequality for quadratic Hecke characters will be an important ingredient of this paper.    The study of the large sieve inequality for characters of a fixed order is of independent interest.  We refer the reader to \cites{DRHB, DRHB1, G&Zhao, G&L, B&Y, BGL}.

\begin{lemma}{\cite[Theorem 1]{Onodera}} \label{quartls}
Suppose that $K= \mq(i)$ or $\mq(\omega)$.  Let $M$, $N$ be positive integers, and let $(a_n)_{n\in \mathbb{N}}$ be an arbitrary sequence of complex numbers, where $n$ runs over $\mz[i]$. Then we have
\begin{equation*}
%%\label{eq:quartic}
 \sumstar_{\substack{m \in \mathcal{O}_K  \\N(m) \leq M}} \left| \ \sumstar_{\substack{n \in \mathcal{O}_K \\N(n) \leq N}} a_n \leg{n}{m} \right|^2
 \ll_{\varepsilon} (M + N)(MN)^{\varepsilon} \sum_{N(n) \leq N} |a_n|^2,
\end{equation*}
   for any $\varepsilon > 0$, where the asterisks indicate that $m$ and $n$ run over square-free elements of $\mathcal{O}_K$ and $(\frac
{\cdot}{m})$ is the quadratic residue symbol.
\end{lemma}

\section{Proof of Theorem~\ref{firstmoment}}

The proofs for $\mq(i)$ and $\mq(\omega)$ are similar.  We will give the full details of the proof for $\mq(i)$ and a sketch of the proof for $\mq(\omega)$.

\subsection{The main term of the first moment}

We have, using \eqref{quadapproxfuneqQi} with $G(s)=1$, that
\begin{align*} 
%%\label{sum1}
   \sumstar_{c \equiv 1 \bmod {16}}L \left( \frac{1}{2},
   \chi_c \right) \Phi\left( \frac{N(c)}{y} \right) &=2\sumstar_{c \equiv 1 \bmod {16}} \ \sum_{0 \neq \mathcal{A} \subset
  O_K} \frac{\chi_c(\mathcal{A})}{N(\mathcal{A})^{1/2}}V \left(\frac{\pi N(\mathcal{A}) }{N(c)^{1/2}} \right)\Phi\left( \frac{N(c)}{y} \right).
\end{align*}

   Since any integral non-zero ideal $\mathcal{A}$ in $\mz[i]$ has a unique generator
$(1+i)^ra$, with $r \in \intz, r \geq 0 , a \in \intz[i], a \equiv 1 \pmod
{(1+i)^3}$, it follows from the quadratic reciprocity law and \eqref{2.5} that
$\chi_{c}(\mathcal{A}) = \chi_a(c)$ (recall our convention that $\chi_1$ is the principal character $\pmod 1$). This allows us to recast the last expression above as
\[  M= 2 \sum_{\substack{r \geq 0 \\ a \equiv 1 \bmod {(1+i)^3}}} \frac{1}{2^{r/2}N(a)^{1/2}}M(r,a),
\; \; \mbox{where} \; \; M(r,a)= \sumstar_{c \equiv 1 \bmod {16}}\chi_a(c)V\left( \frac{\pi 2^rN(a)}{y^{1/2}} \frac{y^{1/2}}{N(c)^{1/2}} \right ) \Phi\left( \frac{N(c)}{y} \right). \]
  Now we use M\"obius inversion (writing $\mu_{[i]}$ for the M\"obius function on $\mz[i]$) to detect the
condition that $c$ is square-free, getting
\[   M(r,a)= \sum_{l \equiv 1 \bmod {(1+i)^3}}\mu_{[i]}(l)\chi_a(l^2) M(l,r,a), \; \; \mbox{with} \; \;   M(l,r,a)= \sum_{c \equiv 1 \bmod {16}}\chi_a(c)V \left( \frac{\pi 2^{r}N(a)}{y^{1/2}} \frac{y^{1/2}}{N(cl^2)^{1/2}} \right )\Phi\left( \frac{N(cl^2)}{y} \right). \]

    By Mellin inversion, we have
\[  V \left( \frac{\pi 2^rN(a)}{y^{1/2}} \frac{y^{1/2}}{N(cl^2)^{1/2}} \right ) \Phi\left( \frac{N(cl^2)}{y} \right) = \frac 1{2\pi i}\int\limits_{(2)} \left( \frac{y}{N(cl^2)} \right)^s \tilde{f}(s) \ \dif s, \; \; \mbox{where} \; \;  \tilde{f}(s)=\int\limits^{\infty}_{0}V \left( \frac{\pi 2^rN(a)}{(xy)^{1/2}} \right) \Phi(x) x^{s-1} \dif x.
\]

Integration by parts and using \eqref{2.07} shows $\tilde{f}(s)$ is a function satisfying the bound for all $\Re(s) > 0$, and $E>0$,
\begin{align}
\label{3.1}
  \tilde{f}(s) \ll (1+|s|)^{-E} \left( 1+\frac{2^rN(a)}{y^{1/2}} \right)^{-E}.
\end{align}

    With this notation, we have
\begin{align*}
  M(l,r,a) = \frac 1{2\pi i}\int\limits_{(2)}\tilde{f}(s) \left( \frac{y}{N(l^2)} \right)^s\sum_{c \equiv 1 \bmod {16}}\frac {\chi_a(c)}{N(c)^s} \dif s.
\end{align*}

   We now use the ray class characters to detect the
condition that $c \equiv 1 \bmod {16}$, getting
\begin{align*}
  M(l,r,a)=\frac {1}{\#h_{(16)}}\sum_{\psi \bmod {16}}\frac {1}{2\pi
   i}\int\limits\limits_{(2)}\tilde{f}(s) \left( \frac{y}{N(l^2)} \right)^s L(s, \psi\chi_a)\dif s,
\end{align*}
   where $\psi$ runs over all ray class characters $\pmod {16}$, $\#h_{(16)}=32$ and
\begin{align*}
   L(s, \psi\chi_a)=\sum_{\mathcal{A} \neq 0} \frac{\psi(\mathcal{A})\chi_a(\mathcal{A})}{N(\mathcal{A})^s}.
\end{align*}

   We estimate $M$ by shifting the contour to the half line. When $\psi\chi_a$ is principal, the Hecke $L$-function
has a pole at $s = 1$. We set $M_0$ to be the contribution to $M$ of these residues, and $M_1$ to be the
remainder. We shall determine $M_0$ first. \newline

    Note that $\psi\chi_a$ is principal if and only if both $\psi$ and $\chi_a$ are principal. Hence $a$ must be a square. We denote $\psi_0$ for the principal
    ray class character $\pmod {16}$. Then we have
\begin{align*}
   L(s, \psi_0\chi_{a^2})=\zeta_{\mq(i)}(s)\prod_{\pi |2a} \left(1-N(\pi)^{-s} \right).
\end{align*}

   Let $c_0=\pi/4$, the residue of $\zeta_{\mq(i)}(s)$ at $s=1$.  Then we have
\begin{align*}
  M_0 &=\frac {2y}{\#h_{(16)}} \sum_{\substack{r \geq 0 \\ a \equiv 1 \bmod {(1+i)^3}}} \frac{1}{2^{r/2}N(a)}\tilde{f}(1)\text{Res}_{s=1}L(s, \psi_0\chi_{a^2})\sum_{l \equiv 1 \bmod (1+i)^3} \frac{\mu_{[i]}(l)\chi_{a^2}(l^2)}{N(l^2)}  \\
  &=\frac {2c_0y}{\#h_{(16)}\zeta_{\mq(i)}(2)} \sum_{\substack{r \geq 0 \\ a \equiv 1 \bmod {(1+i)^3}}} \frac{1}{2^{r/2}N(a)}\tilde{f}(1)\prod_{\pi |
  2a}\left( 1-N(\pi)^{-1} \right)\prod_{\pi | (2a)} \left( 1-N(\pi)^{-2} \right)^{-1} \\
  &=\frac {2c_0y}{\#h_{(16)}\zeta_{\mq(i)}(2)} \sum_{\substack{r \geq 0 \\ a \equiv 1 \bmod {(1+i)^3}}} \frac{1}{2^{r/2}N(a)}\tilde{f}(1)\prod_{\pi |
  2a} \left( 1+N(\pi)^{-1} \right)^{-1}.
\end{align*}

   Note that
\begin{align*}
  \prod_{\pi | 2a}(1+N(\pi)^{-1})^{-1}=\sum_{(d) | (2a)}\frac {\mu_{[i]}(d)}{\sigma(d)},
\end{align*}
   where $\sigma(d)$ denotes the sum of the norms of the integral ideal divisors of $(d)$. \newline

   Applying this, we see that
\begin{align*}
  \sum_{\substack{N(a) \leq x \\ a \equiv 1 \bmod {(1+i)^3}}} \frac{1}{N(a)}\prod_{\pi |
  2a} \left( 1+N(\pi)^{-1} \right)^{-1}=\sum_{\substack{ (d) \\ N(d) \leq 2x}}\frac {\mu_{[i]}(d)}{\sigma(d)}\sum_{\substack{a \equiv 1 \bmod {(1+i)^3} \\ \frac {d}{(2,d)}|a,\  N(a) \leq x}} \frac{1}{N(a)}.
\end{align*}

    Note the following result from the Gauss circle problem,
\begin{align*}
  \sum_{\substack{N(a) \leq x \\ a \equiv 1 \bmod {(1+i)^3}}} 1 =\frac {\pi}{8}x+O(x^{\theta}).
\end{align*}
    Here one can take $\theta$ to be $131/416$ (see \cite{Huxley1}). \newline

    Applying this and partial summation, we get
\begin{align*}
  \sum_{\substack{N(a) \leq x \\ a \equiv 1 \bmod {(1+i)^3}}} \frac{1}{N(a)} = \frac {\pi}{8}\log x+C_0+O(x^{\theta-1}),
\end{align*}
    where $C_0$ is a constant. \newline

    It follows that
\begin{align*}
  \sum_{\substack{N(a) \leq x \\ a \equiv 1 \bmod {(1+i)^3}}} \frac{1}{N(a)}\prod_{\pi |
  2a} \left( 1+N(\pi)^{-1} \right)^{-1}=\frac {\pi}{12}A_{\mq(i)} \log x+C_1+O(x^{\theta-1}),
\end{align*}
   where $A_{\mq(i)}$, $C_1$ are constants with $A_{\mq(i)}$ defined in \eqref{1.4}. \newline

   We then deduce that
\[ \sum_{a \equiv 1 \bmod {(1+i)^3}} \frac{1}{N(a)}\tilde{f}(1)\prod_{\pi |
  2a} \left( 1+N(\pi)^{-1} \right)^{-1} = \int\limits^{2}_{1}\Phi(x) \sum_{a \equiv 1 \bmod {(1+i)^3}} \frac{1}{N(a)}\prod_{\pi |
  2a}(1+N(\pi)^{-1})^{-1}V \left( \frac{\pi 2^rN(a)^2}{(xy)^{1/2}} \right) \dif x. \]

Applying partial summation and \eqref{2.07}, we get
\begin{align*}
  & \sum_{a \equiv 1 \bmod {(1+i)^3}} \frac{1}{N(a)}\prod_{\pi |
  2a}(1+N(\pi)^{-1})^{-1}V \left( \frac{\pi 2^rN(a)^2}{(xy)^{1/2}} \right)\\
 & =\begin{cases}
  \displaystyle \frac {\pi}{12}A_{\mq(i)}\log \frac {(xy)^{1/4}}{\pi^{1/2} 2^{r/2}}+C_2+O\left( \left( \frac {\pi^{1/2} 2^{r/2}}{(xy)^{1/4}} \right)^{1-\theta} \right) \qquad & (xy)^{1/2} > \pi 2^{r/2} , \\ \\
\displaystyle     O\left( \left( \frac {\pi^{1/2} 2^{r/2}}{(xy)^{1/4}} \right)^{-1} \right)  \qquad & \displaystyle (xy)^{1/2} \leq \pi 2^{r/2}  ,
    \end{cases}
\end{align*}
   with some constant $C_2$. \newline

   We then conclude that by a straightforward calculation (we may assume that $y$ is large),
\begin{align} \label{M0formula}
  M_0 =\frac {(2+\sqrt{2})c_0\pi A_{\mq(i)}}{24 \#h_{(16)}\zeta_{\mq(i)}(2)}\hat{\Phi}(0) y\log y +C_{\mq(i)} \hat{\Phi}(0)y+O(y^{(3+\theta)/4}),
\end{align}
   where $C_{\mq(i)}$ is the same constant $C_{\mq(i)}$ appearing in \eqref{1stmom}.

\subsection{The remainder terms of the first moment}

   To treat $M_1$, we bound everything by absolute values and use \eqref{3.1} to get that for any $E>0$,
\begin{align}
\label{3.2}
   M_1 \ll y^{1/2} \sum_{N(l) \ll \sqrt{y}}\frac {1}{N(l)} \sum_{\psi \bmod {16}} \sum_{\substack{r \geq 0 \\ a \equiv 1 \bmod {(1+i)^3}}} N(a)^{-1/2} \left( 1+\frac{2^rN(a)}{y^{1/2}} \right)^{-E} \int\limits^{\infty}_{-\infty} \left| L\left(\frac 12+it, \psi\chi_a \right) \right| (1+|t|)^{-E} \dif t.
\end{align}
  Now it follows easily from the Cauchy-Schwarz inequality and \eqref{3.4} that
\begin{align*}
%%\label{3.04}
    \sum_{\substack{a \equiv 1 \bmod {(1+i)^3} \\ N(a) \leq N}}  N(a)^{-1/2} \left| L\left( \frac 12+it, \psi\chi_a \right) \right| \ll (N(1+|t|))^{1/2+\epsilon}.
\end{align*}

  Applying this in \eqref{3.2} and note that we can restrict the sum over $r, a$ to be $2^r N(a) \leq y^{1/2+\epsilon}$, we immediately deduce that
\begin{align*}
   M_1 \ll y^{3/4+\epsilon}.
\end{align*}

Combining the results for $M_0$ and $M_1$, we obtain the result of Theorem \ref{firstmoment} for $\mq(i)$.

\subsection{The proof for $\mq(\omega)$}  As stated before, the proof for $\mq(\omega)$ is very similar.  Starting with the approximate functional equation, \eqref{quadapproxfuneqQomega}, with $G(s)=1$, we have
\begin{align*}
%%\label{sum1}
   M=\sumstar_{c \equiv 1 \bmod {36}}L \left( \frac{1}{2},
   \chi_c \right) \Phi \left( \frac{N(c)}{y} \right) &=2\sumstar_{c \equiv 1 \bmod {36}} \ \sum_{0 \neq \mathcal{A} \subset
  O_K} \frac{\chi_c(\mathcal{A})}{N(\mathcal{A})^{1/2}}V \left(\frac{2\pi N(\mathcal{A}) }{(3N(c))^{1/2}} \right) \Phi \left( \frac{N(c)}{y} \right).
\end{align*}

Since any integral non-zero ideal $\mathcal{A}$ in $\mz[\omega]$ has a unique generator
$2^{r_1}(1-\omega)^{r_2}a$, with $r_1, r_2\in \intz, r_1,r_2 \geq 0 , a \in \intz[\omega]$,  $(a, 2)=1, a \equiv 1 \pmod 3$, it follows from \eqref{2.5} and the definition of $\chi^{(a)}$ that $\chi_{c}(\mathcal{A}) = \chi_c(a)$. \newline

The above discussions allow us to recast $M$ as
\[  M= 2 \sum_{\substack{r_1,r_2 \geq 0 \\ (a,2)=1 \\ a \equiv 1 \bmod 3}} \frac{1}{2^{r_1}3^{r_2/2}N(a)^{1/2}}M(r,a),\]
where
\[ M(r,a)= \sumstar_{c \equiv 1 \bmod {36}}\chi^{(a)}(c)V\left( \frac{\pi 2^{2r_1+1}3^{r_2-1/2}N(a)}{y^{1/2}} \frac{y^{1/2}}{N(c)^{1/2}} \right ) \Phi \left( \frac{N(c)}{y} \right). \]
Now the proof goes along the same arguments as those in the proof for $\mq(i)$.  The main term, the analogue of \eqref{M0formula}, that emerges from this computation is
\[   \frac {(1+\sqrt{3})\pi A_{\mq(\omega)} }{20 \#h_{(36)}\zeta_{\mq(\omega)}(2)} \mathrm{Res}_{s=1} \zeta_{\mq(\omega)}(s)   \widehat{\Phi}(0) y\log y +C_{\mq(\omega)} \widehat{\Phi}(0)y+O(y^{(3+\theta)/4}), \]
   where $C_{\mq(\omega)}$ is the same constant $C_{\mq(\omega)}$ appearing in \eqref{1stmom}, 
\[  \mathrm{Res}_{s=1} \zeta_{\mq(\omega)}(s)  = \frac{\sqrt{3}}{9} \pi , \quad   \#h_{(36)} = 108. \]   
This completes the proof of the theorem.

\section{Proof of Theorem~\ref{secmom}}

The proofs for $\mq(i)$ and $\mq(\omega)$ are very similar and we only give the details for $\mq(i)$ here. \newline

Let $a \equiv 1 \pmod {(1+i)^3}$ and $\psi$ be a ray class character $\pmod {16}$. Let $\chi_a$ be the Hecke character $\pmod {16a}$ of trivial infinite type. For any ideal $(c)$ co-prime to $(1+i)$, with $c$ being the unique generator of $(c)$ satisfying $c \equiv 1 \pmod {(1+i)^3}$, $\chi_a((c))$ is defined as $\chi_a((c))=\leg {a}{c}$ (see \cite[Example 2, p. 62]{HIEK}). Now we show that
\begin{align}
\label{3.4}
    \sum_{\substack{a \equiv 1 \bmod {(1+i)^3} \\ N(a) \leq N}} \left| L \left( \frac 12+it, \psi\chi_a \right) \right|^2 \ll \left( N(1+|t|) \right)^{1+\epsilon}.
\end{align}

    As the proof is similar to the proof of \cite[Theorem 1.3]{B&Y} and \cite[Corollary 1.4]{BGL}, we only sketch the arguments.
    Write $a=a_1a^2_2$ with $a_1, a_2 \equiv 1 \pmod {(1+i)^3}$ and $a_1$ square-free. Then $\chi_a$ equals $\chi_{a_1}$ multiplied by a principal character whose conductor divides $a_2$.  We may further assume that $\psi\chi_{a_1}$ is primitive. Then \eqref{approxfuneq} is valid with $G(s)=e^{-s^2}$. By inserting \eqref{approxfuneq} into the left-hand side expression in \eqref{3.4} with $x=2(N(m_{a_1}))^{1/2}$, where $m_{a_1}$ is the conductor of $\psi\chi_{a_1}$, and applying the Cauchy-Schwarz inequality, we see that it suffices to bound
\begin{align*}
   \sum_{\substack{a_2 \equiv 1 \bmod {(1+i)^3} \\ N(a_2)^2 \leq N}} \ \sumstar_{\substack{a_1 \equiv 1 \bmod {(1+i)^3} \\ N(a_1) \leq N/N(a_2)^2}} \left| \sum_{0 \neq \mathcal{A} \subset
  \mathcal{O}_K}\frac{\psi\chi_{a_1}(\mathcal{A})}{N(\mathcal{A})^{1/2+it}}V_t\left (\frac{2\pi N(\mathcal{A})}{N(m_{a_1})^{1/2}}  \right) \right|^2,
\end{align*}
     where $\sum^*$ denotes summation over square-free elements of $\mz[i]$. \newline
     
     In view of \eqref{2.15}, we may truncate the sum over $\mathcal{A}$ above to $N(\mathcal{A}) \leq (N^{1/2}(1+|t|))^{1+\epsilon}$. By shifting the contour in \eqref{2.14} to the $\epsilon$ line and write $s=\epsilon+iw$, it suffices to bound
\begin{align*}
%%\label{3.9}
    \sum_{\substack{a_2 \equiv 1 \bmod {(1+i)^3} \\ N(a_2)^2 \leq N}}\sumstar_{\substack{a_1 \equiv 1 \bmod {(1+i)^3} \\ N(a_1) \leq N/N(a_2)^2}} \left| \sum_{\substack {0 \neq \mathcal{A} \subset
  \mathcal{O}_K \\ N(\mathcal{A}) \leq (N^{1/2}(1+|t|))^{1+\epsilon}}}\frac{\psi\chi_{a_1}(\mathcal{A})}{N(\mathcal{A})^{1/2+\epsilon+it+iw}} \right|^2.
\end{align*}

    In the inner sum above, writing $\mathcal{A}=(1+i)^r\mathcal{A}_1\mathcal{A}^2_2$ with $\mathcal{A}_1$, $\mathcal{A}_2$ co-prime to $1+i$, $\mathcal{A}_1$ square-free and using the Cauchy-Schwarz inequality, it is enough to estimate
\begin{align*}
   \sum_{\substack{a_2 \equiv 1 \bmod {(1+i)^3} \\ N(a_2)^2 \leq N}}  \sum_{\substack {r \geq 0\\ 2^r \leq (N^{1/2}(1+|t|)^{1/2})^{1+\epsilon}}} & \frac {1}{2^{r(1+2\epsilon)/2}} \sum_{\substack {0 \neq \mathcal{A}_2 \subset
  \mathcal{O}_K \\ N(\mathcal{A}_2)^2 \leq (N^{1/2}(1+|t|))^{1+\epsilon}/2^r}}\frac{1}{N(\mathcal{A}_2)^{1+2\epsilon}} \\
  &\times \sumstar_{\substack{a_1 \equiv 1 \bmod {(1+i)^3} \\ N(a_1) \leq N/N(a_2)^2}} \left| \ \sumstar_{\substack {0 \neq \mathcal{A}_1 \subset
  \mathcal{O}_K \\ N(\mathcal{A}_1) \leq (N^{1/2}(1+|t|))^{1+\epsilon}/(N(\mathcal{A}^2_2)2^r)}}\frac{\psi\chi_{a_1}(\mathcal{A}_1)}{N(\mathcal{A}_1)^{1/2+\epsilon+it+iw}} \right|^2.
\end{align*}

   We now apply the large sieve inequality and arrive at
\begin{align*}
  & \sumstar_{\substack{a_1 \equiv 1 \bmod {(1+i)^3} \\ N(a_1) \leq N/N(a_2)^2}} \left| \ \sumstar_{\substack {0 \neq \mathcal{A}_1 \subset
  \mathcal{O}_K \\ N(\mathcal{A}_1) \leq (N^{1/2}(1+|t|))^{1+\epsilon}/(N(\mathcal{A}^2_2)2^r)}}\frac{\psi\chi_{a_1}(\mathcal{A}_1)}{N(\mathcal{A}_1)^{1/2+\epsilon+it+iw}}\right|^2 \ll N^{\epsilon} \left( \frac {N}{N(a_2)^2}+ \frac {(N^{1/2}(1+|t|))^{1+\epsilon}}{(N(\mathcal{A}^2_2)2^r)} \right).
\end{align*}

    As the sums over $a_2, r$ and $\mathcal{A}_2$ all converge, we conclude that
\begin{align*}
   \sum_{\substack{a_2 \equiv 1 \bmod {(1+i)^3} \\ N(a_2)^2 \leq N}}\sumstar_{\substack{a_1 \equiv 1 \bmod {(1+i)^3} \\ N(a_1) \leq N/N(a_2)^2}} \left| \sum_{\substack {0 \neq \mathcal{A} \subset
  \mathcal{O}_K \\ N(\mathcal{A}) \leq (N^{1/2}(1+|t|))^{1+\epsilon}}}\frac{\psi\chi_{a_1}(\mathcal{A})}{N(\mathcal{A})^{1/2+\epsilon+it+iw}} \right|^2
   \ll \left( N(1+|t|) \right)^{1+\epsilon}.
\end{align*}

     This establishes \eqref{3.4} from which Theorem~\ref{secmom} follows readily. \newline

\noindent{\bf Acknowledgments.} P. G. is supported in part by NSFC grant 11871082 and L. Z. by the FRG grant PS43707.  Parts of this work were done when P. G. visited the University of New South Wales (UNSW) in June 2017. He wishes to thank UNSW for the invitation, financial support and warm hospitality during his pleasant stay.

\bibliography{biblio}
\bibliographystyle{amsxport}

\vspace*{.5cm}

\noindent\begin{tabular}{p{8cm}p{8cm}}
School of Mathematics and Systems Science & School of Mathematics and Statistics \\
Beihang University & University of New South Wales \\
Beijing 100191 China & Sydney NSW 2052 Australia \\
Email: {\tt penggao@buaa.edu.cn} & Email: {\tt l.zhao@unsw.edu.au} \\
\end{tabular}

\end{document}